# SPARSISTENCY AND RATES OF CONVERGENCE IN LARGE COVARIANCE MATRIX ESTIMATION[1]

BY CLIFFORD LAM AND JIANQING FAN

*London School of Economics and Political Science and Princeton University*

This paper studies the sparsistency and rates of convergence for estimating sparse covariance and precision matrices based on penalized likelihood with nonconvex penalty functions. Here, sparsistency refers to the property that all parameters that are zero are actually estimated as zero with probability tending to one. Depending on the case of applications, sparsity priori may occur on the covariance matrix, its inverse or its Cholesky decomposition. We study these three sparsity exploration problems under a unified framework with a general penalty function. We show that the rates of convergence for these problems under the Frobenius norm are of order $(s_n \log p_n/n)^{1/2}$, where $s_n$ is the number of nonzero elements, $p_n$ is the size of the covariance matrix and $n$ is the sample size. This explicitly spells out the contribution of high-dimensionality is merely of a logarithmic factor. The conditions on the rate with which the tuning parameter $\lambda_n$ goes to 0 have been made explicit and compared under different penalties. As a result, for the $L_1$-penalty, to guarantee the sparsistency and optimal rate of convergence, the number of nonzero elements should be small: $s'_n = O(p_n)$ at most, among $O(p_n^2)$ parameters, for estimating sparse covariance or correlation matrix, sparse precision or inverse correlation matrix or sparse Cholesky factor, where $s'_n$ is the number of the nonzero elements on the off-diagonal entries. On the other hand, using the SCAD or hard-thresholding penalty functions, there is no such a restriction.

**1. Introduction.** Covariance matrix estimation is a common statistical problem in many scientific applications. For example, in financial risk assessment or longitudinal study, an input of covariance matrix $\mathbf{\Sigma}$ is needed,

Received February 2009; revised March 2009.
[1]Supported by NSF Grants DMS-03-54223, DMS-07-04337 and NIH Grant R01-GM072611.

whereas an inverse of the covariance matrix, the precision matrix $\mathbf{\Sigma}^{-1}$, is required for optimal portfolio selection, linear discriminant analysis or graphical network models. Yet, the number of parameters in the covariance matrix grows quickly with dimensionality. Depending on the applications, the sparsity of the covariance matrix or precision matrix is frequently imposed to strike a balance between biases and variances. For example, in longitudinal data analysis [see, e.g., Diggle and Verbyla (1998), or Bickel and Levina (2008b)], it is reasonable to assume that remote data in time are weakly correlated, whereas in Gaussian graphical models, the sparsity of the precision matrix is a reasonable assumption [Dempster (1972)].

This initiates a series of researches focusing on the parsimony of a covariance matrix. Smith and Kohn (2002) used priors which admit zeros on the off-diagonal elements of the Cholesky factor of the precision matrix $\mathbf{\Omega} = \mathbf{\Sigma}^{-1}$, while Wong, Carter and Kohn (2003) used zero-admitting prior directly on the off-diagonal elements of $\mathbf{\Omega}$ to achieve parsimony. Wu and Pourahmadi (2003) used the Modified Cholesky Decomposition (MCD) to find a banded structure for $\mathbf{\Omega}$ nonparametrically for longitudinal data. Bickel and Levina (2008b) developed consistency theories on banding methods for longitudinal data, for both $\mathbf{\Sigma}$ and $\mathbf{\Omega}$.

Various authors have used penalized likelihood methods to achieve parsimony on covariance selection. Fan and Peng (2004) has laid down a general framework for penalized likelihood with diverging dimensionality, with general conditions for the oracle property stated and proved. However, it is not clear whether it is applicable to the specific case of covariance matrix estimation. In particular, they did not link the dimensionality $p_n$ with the number of nonzero elements $s_n$ in the true covariance matrix $\mathbf{\Sigma}_0$, or the precision matrix $\mathbf{\Omega}_0$. A direct application of their results to our setting can only handle a relatively small covariance matrix of size $p_n = o(n^{1/10})$.

Recently, there is a surge of interest on the estimation of sparse covariance matrix or precision matrix using penalized likelihood method. Huang et al. (2006) used the LASSO on the off-diagonal elements of the Cholesky factor from MCD, while Meinshausen and Bühlmann (2006), d'Aspremont, Banerjee and El Ghaoui (2008) and Yuan and Lin (2007) used different LASSO algorithms to select zero elements in the precision matrix. A novel penalty called the nested LASSO was constructed in Levina, Rothman and Zhu (2008) to penalize off-diagonal elements. Thresholding the sample covariance matrix in high-dimensional setting was thoroughly studied by El Karoui (2008) and Bickel and Levina (2008a) and Cai, Zhang and Zhou (2008) with remarkable results for high-dimensional applications. However, it is not directly applicable to estimating sparse precision matrix when the dimensionality $p_n$ is greater than the sample size $n$. Wagaman and Levina (2008) proposed an Isomap method for discovering meaningful orderings of variables based on their correlations that result in block-diagonal or banded correlation structure, resulting in an Isoband estimator. A permutation invariant estimator, called SPICE, was proposed in Rothman et al. (2008) based on penalized likelihood with $L_1$-penalty on the off-diagonal elements for the precision matrix. They obtained remarkable results on the rates of convergence. The rate for estimating $\mathbf{\Omega}$ under the Frobenius norm is of order $(s_n \log p_n/n)^{1/2}$, with dimensionality cost only a logarithmic factor in the overall mean-square error, where $s_n = p_n + s_{n1}$, $p_n$ is the number of the diagonal elements and $s_{n1}$ is the number of the nonzero off-diagonal entries. However, such a rate of convergence neither addresses explicitly the issues of sparsistency such as those in Fan and Li (2001) and Zhao and Yu (2006), nor the bias issues due to the $L_1$-penalty and the sampling distribution of the estimated nonzero elements. These are the core issues of the study. By sparsistency, we mean the property that all parameters that are zero are actually estimated as zero with probability tending to one, a weaker requirement than that of Ravikumar et al. (2007).

In this paper, we investigate the aforementioned problems using the penalized pseudo-likelihood method. Assume a random sample $\{\mathbf{y}_i\}_{1 \le i \le n}$ with mean zero and covariance matrix $\mathbf{\Sigma}_0$, satisfying some sub-Gaussian tails conditions as specified in Lemma 2 (see Section 5). The sparsity of the true precision matrix $\mathbf{\Omega}_0$ can be explored by maximizing the Gaussian quasi-likelihood or equivalently minimizing

$$(1.1) \qquad q_1(\mathbf{\Omega}) = \text{tr}(\mathbf{S}\mathbf{\Omega}) - \log|\mathbf{\Omega}| + \sum_{i \ne j} p_{\lambda_{n1}}(|\omega_{ij}|),$$

which is the penalized negative log-likelihood if the data is Gaussian. The matrix $\mathbf{S} = n^{-1} \sum_{i=1}^n \mathbf{y}_i \mathbf{y}_i^T$ is the sample covariance matrix, $\mathbf{\Omega} = (\omega_{ij})$, and $p_{\lambda_{n1}}(\cdot)$ is a penalty function, depending on a regularization parameter $\lambda_{n1}$, which can be nonconvex. For instance, the $L_1$-penalty $p_\lambda(\theta) = \lambda|\theta|$ is convex, while the hard-thresholding penalty defined by $p_\lambda(\theta) = \lambda^2 - (|\theta| - \lambda)^2 \mathbf{1}_{\{|\theta| < \lambda\}}$, and the SCAD penalty defined by

$$(1.2) \quad p'_\lambda(\theta) = \lambda \mathbf{1}_{\{\theta \le \lambda\}} + (a\lambda - \theta)_+ \mathbf{1}_{\{\theta > \lambda\}}/(a-1) \qquad \text{for some } a > 2,$$

are folded-concave. Nonconvex penalty is introduced to reduce bias when the true parameter has a relatively large magnitude. For example, the SCAD penalty remains constant when $\theta$ is large, while the $L_1$-penalty grows linearly with $\theta$. See Fan and Li (2001) for a detailed account of this and other advantages of such a penalty function. The computation can be done via the local linear approximation [Zhou and Li (2008) and Fan, Feng and Wu (2009)]; see Section 2.1 for additional details.

Similarly, the sparsity of the true covariance matrix $\mathbf{\Sigma}_0$ can be explored by minimizing

$$(1.3) \qquad q_2(\mathbf{\Sigma}) = \text{tr}(\mathbf{S}\mathbf{\Sigma}^{-1}) + \log|\mathbf{\Sigma}| + \sum_{i \ne j} p_{\lambda_{n2}}(|\sigma_{ij}|),$$

where $\boldsymbol{\Sigma} = (\sigma_{ij})$. Note that we only penalize the off-diagonal elements of $\boldsymbol{\Sigma}$ or $\boldsymbol{\Omega}$ in the aforementioned two methods, since the diagonal elements of $\boldsymbol{\Sigma}_0$ and $\boldsymbol{\Omega}_0$ do not vanish.

In studying a sparse covariance or precision matrix, it is important to distinguish between the diagonal and off-diagonal elements, since the diagonal elements are always positive and they contribute to the overall mean-squares errors. For example, the true correlation matrix, denoted by $\boldsymbol{\Gamma}_0$, has the same sparsity structure as $\boldsymbol{\Sigma}_0$ without the need to estimating its diagonal elements. In view of this fact, we introduce a revised method (3.2) to take this advantage. It turns out that the correlation matrix can be estimated with a faster rate of convergence, at $(s_{n1} \log p_n/n)^{1/2}$ instead of $((p_n + s_{n1}) \log p_n/n)^{1/2}$, where $s_{n1}$ is the number of nonzero correlation coefficients. We can take similar advantages over the estimation of the true inverse correlation matrix, denoted by $\boldsymbol{\Psi}_0$. See Section 2.5. This is an extension of the work of Rothman et al. (2008) using the $L_1$-penalty. Such an extension is important since the nonconcave penalized likelihood ameliorates the bias problem for the $L_1$-penalized likelihood.

The bias issues of the commonly used $L_1$-penalty, or LASSO, can be seen from our theoretical results. In fact, due to the bias of LASSO, an upper bounded of $\lambda_{ni}$ is needed in order to achieve fast rate of convergence. On the other hand, a lower bound is required in order to achieve sparsity of estimated precision or covariance matrices. This is in fact one of the motivations for introducing nonconvex penalty functions in Fan and Li (2001) and Fan and Peng (2004), but we state and prove the explicit rates in the current context. In particular, we demonstrate that the $L_1$-penalized estimator can achieve simultaneously the optimal rate of convergence and sparsistency for estimation of $\boldsymbol{\Sigma}_0$ or $\boldsymbol{\Omega}_0$ when the number of nonzero elements in the off-diagonal entries are no larger than $O(p_n)$, but not guaranteed so otherwise. On the other hand, using the nonconvex penalties like the SCAD or hard-thresholding penalty, such an extra restriction is not needed.

We also compare two different formulations of penalized likelihood using the modified Cholesky decomposition, exploring their respective rates of convergence and sparsity properties.

Throughout this paper, we use $\lambda_{\min}(A)$, $\lambda_{\max}(A)$ and $\operatorname{tr}(A)$ to denote the minimum eigenvalue, maximum eigenvalue, and trace of a symmetric matrix $A$, respectively. For a matrix $B$, we define the operator norm and the Frobenius norm, respectively, as $\|B\| = \lambda_{\max}^{1/2}(B^T B)$ and $\|B\|_F = \operatorname{tr}^{1/2}(B^T B)$.

**2. Estimation of sparse precision matrix.** In this section, we present the analysis of (1.1) for estimating a sparse precision matrix. Before this, let us first present an algorithm for computing the nonconcave maximum (pseudo)-likelihood estimator and then state the conditions needed for our technical results.

2.1. *Algorithm based on iterated reweighted $L_1$-penalty.* The computation of the nonconcave maximum likelihood problems can be solved by a sequence of $L_1$-penalized likelihood problems via local linear approximation [Zou and Li (2008) and Fan, Feng and Wu (2009)]. For example, given the current estimate $\Omega_k = (\omega_{ij,k})$, by the local linear approximation to the penalty function,

$$
\begin{aligned}
q_1(\boldsymbol{\Omega}) &\approx \operatorname{tr}(\mathbf{S}\boldsymbol{\Omega}) - \log|\boldsymbol{\Omega}| \\
&\quad + \sum_{i \neq j}[p_{\lambda_{n1}}(|\omega_{ij,k}|) + p'_{\lambda_{n1}}(|\omega_{ij,k}|)(|\omega_{ij}| - |\omega_{ij,k}|)].
\end{aligned}
\tag{2.1}
$$

Hence, $\boldsymbol{\Omega}_{k+1}$ should be taken to maximize the right-hand side of (2.1):

$$
\boldsymbol{\Omega}_{k+1} = \arg\max_{\boldsymbol{\Omega}}\left[\operatorname{tr}(\mathbf{S}\boldsymbol{\Omega}) - \log|\boldsymbol{\Omega}| + \sum_{i \neq j} p'_{\lambda_{n1}}(|\omega_{ij,k}|)|\omega_{ij}|\right],
\tag{2.2}
$$

after ignoring the two constant terms. Problem (2.2) is the weighted penalized $L_1$-likelihood. In particular, if we take the most primitive initial value $\boldsymbol{\Omega}_0 = \mathbf{0}$, then

$$
\boldsymbol{\Omega}_1 = \arg\max_{\boldsymbol{\Omega}}\left[\operatorname{tr}(\mathbf{S}\boldsymbol{\Omega}) - \log|\boldsymbol{\Omega}| + \lambda_{n1} \sum_{i \neq j} |\omega_{ij}|\right],
$$

is already a good estimator. Iterations of (2.2) reduces the biases of the estimator, as larger estimated coefficients in the previous iterations receive less penalty. In fact, in a different setup, Zou and Li (2008) showed that one iteration of such a procedure is sufficient as long as the initial values are good enough.

Fan, Feng and Wu (2009) has implemented the above algorithm for optimizing (1.1). They have also demonstrated in Section 2.2 in their paper how to utilize the graphical lasso algorithm of Friedman, Hastie and Tibshirani (2008), which is essentially a group coordinate descent procedure, to solve problem (2.2) quickly, even when $p_n > n$. Such a group coordinate decent algorithm was also used by Meier, van de Geer and Bühlmann (2008) to solve the group LASSO problem. Thus iteratively, (2.2), and hence (1.1), can be solved quickly with the graphical lasso algorithm. See also Zhang (2007) for a general solution to the folded-concave penalized least-squares problem. The following is a brief summary of the numerical results in Fan, Feng and Wu (2009).

2.2. *Some numerical results.* We give a brief summary of a breast cancer data analysis with $p_n > n$ considered in Fan, Feng and Wu (2009). For full details, please refer to Section 3.2 of Fan, Feng and Wu (2009). Other simulation results are also in Section 4 in their paper.

*Breast cancer data.* Normalized gene expression data from 130 patients with stage I–III breast cancers are analyzed, with 33 of them belong to class 1 and 97 belong to class 2. The aim is to assess prediction accuracy in predicting which class a patient will belong to, using a set of pre-selected genes ($p_n = 110$, chosen by $t$-tests) as gene expression profile data. The data is randomly divided into training ($n = 109$) and testing sets. The mean vector for the genes expression levels is obtained from the training data, as well as the associated inverse covariance matrix estimated using LASSO, adaptive LASSO and SCAD penalties as three different regularization methods. A linear discriminant score is then calculated for each regularization method and applied to the testing set to predict if a patient belongs to class 1 or 2. This is repeated 100 times.

On average, the estimated precision matrix $\hat{\boldsymbol{\Omega}}$ using LASSO has many more nonzeros than that using SCAD (3923 versus 674). This is not surprising when we look at (2.3) in our paper, where the $L_1$-penalty imposes an upper bound on the tuning parameter $\lambda_{n1}$ for consistency, which links to reducing the bias in the estimation. This makes the $\lambda_{n1}$ in practice too small to set many of the elements in $\hat{\boldsymbol{\Omega}}$ to zero. While we do not know which elements in the true $\boldsymbol{\Omega}$ are zero, the large number of nonzero elements in the $L_1$-penalized estimator seems spurious, and the resulting gene network is not easy to interpret.

On the other hand, SCAD-penalized estimator has a much smaller number of nonzero elements, since the tuning parameter $\lambda_{n1}$ is not bounded above under consistency of the resulting estimator. This makes the resulting gene network easier to interpret, with some clusters of genes identified.

Also, classification results on the testing set using the SCAD penalty for precision matrix estimation is better than that using the $L_1$-penalty, in the sense that the specificity (#True Negative/#class 2) is higher (0.794 to 0.768) while the sensitivity (#True Positive/#class 1) is similar to that using $L_1$-penalized precision matrix estimator.

2.3. *Technical conditions.* We now introduce some notation and present regularity conditions for the rate of convergence and sparsistency.

Let $S_1 = \{(i,j) : \omega_{ij}^0 \neq 0\}$, where $\boldsymbol{\Omega}_0 = (\omega_{ij}^0)$ is the true precision matrix. Denote by $s_{n1} = |S_1| - p_n$, which is the number of nonzero elements in the off-diagonal entries of $\boldsymbol{\Omega}_0$. Define

$$a_{n1} = \max_{(i,j) \in S_1} p'_{\lambda_{n1}}(|\omega_{ij}^0|), \qquad b_{n1} = \max_{(i,j) \in S_1} p''_{\lambda_{n1}}(|\omega_{ij}^0|).$$

The term $a_{n1}$ is related to the asymptotic bias of the penalized likelihood estimate due to penalization. Note that for $L_1$-penalty, $a_{n1} = \lambda_{n1}$ and $b_{n1} = 0$, whereas for SCAD, $a_{n1} = b_{n1} = 0$ for sufficiently large $n$ under the last assumption of condition (B) below.

We assume the following regularity conditions:

(A) There are constants $\tau_1$ and $\tau_2$ such that
$$0 < \tau_1 < \lambda_{\min}(\boldsymbol{\Sigma}_0) \leq \lambda_{\max}(\boldsymbol{\Sigma}_0) < \tau_2 < \infty \qquad \text{for all } n.$$

(B) $a_{n1} = O(\{(1+p_n/(s_{n1}+1))\}(\log p_n/n)^{1/2})$, $b_{n1} = o(1)$, and $\min_{(i,j) \in S_1} |\omega_{ij}^0|/\lambda_{n1} \to \infty$ as $n \to \infty$.

(C) The penalty $p_\lambda(\cdot)$ is singular at the origin, with $\lim_{t \downarrow 0} p_\lambda(t)/(\lambda t) = k > 0$.

(D) There are constants $C$ and $D$ such that, when $\theta_1, \theta_2 > C\lambda_{n1}, |p''_{\lambda_{n1}}(\theta_1) - p''_{\lambda_{n1}}(\theta_2)| \leq D|\theta_1 - \theta_2|$.

Condition (A) bounds uniformly the eigenvalues of $\boldsymbol{\Sigma}_0$, which facilitates the proof of consistency. It also includes a wide class of covariance matrices as noted in Bickel and Levina (2008b). The rates $a_{n1}$ and $b_{n1}$ in condition (B) are also needed for proving consistency. If they are too large, the bias due to penalty can dominate the variance from the likelihood, resulting in poor estimates.

The last requirement in condition (B) states the rate at which the nonzero parameters should be distinguished from zero asymptotically. It is not explicitly needed in the proofs, but for asymptotically unbiased penalty functions, this is a necessary condition so that $a_{n1}$ and $b_{n1}$ are converging to zero fast enough as needed in the first part of condition (B). In particular, for the SCAD and hard-thresholding penalty functions, this condition implies that $a_{n1} = b_{n1} = 0$ exactly for sufficiently large $n$, thus allowing a flexible choice of $\lambda_{n1}$. For the SCAD penalty (1.2), the condition can be relaxed as $\min_{(i,j) \in S_1} |\omega_{ij}^0|/\lambda_{n1} > a$.

The singularity in condition (C) gives sparsity in the estimates [Fan and Li (2001)]. Finally, condition (D) is a smoothing condition for the penalty function, and is needed in proving asymptotic normality. The SCAD penalty, for instance, satisfies this condition by choosing the constant $D$, independent of $n$, to be large enough.

2.4. *Properties of sparse precision matrix estimation.* Minimizing (1.1) involves nonconvex minimization, and we need to prove that there exists a local minimizer $\hat{\boldsymbol{\Omega}}$ for the minimization problem with a certain rate of convergence, which is given under the Frobenius norm. The proof is given in Section 5. It is similar to the one given in Rothman et al. (2008), but now the penalty function is nonconvex.

THEOREM 1 (Rate of convergence). *Under regularity conditions* (A)–(D), *if* $(p_n + s_{n1}) \log p_n/n = O(\lambda_{n1}^2)$ *and* $(p_n + s_{n1})(\log p_n)^k/n = O(1)$ *for some* $k > 1$, *then there exists a local minimizer* $\hat{\boldsymbol{\Omega}}$ *such that* $\|\hat{\boldsymbol{\Omega}} - \boldsymbol{\Omega}_0\|_F^2 = O_P\{(p_n + s_{n1}) \log p_n/n\}$. *For the $L_1$-penalty, we only need* $\log p_n/n = O(\lambda_{n1}^2)$.

The proofs of this theorem and others are relegated to Section 5 so that readers can get more quickly what the results are. As in Fan and Li (2001), the asymptotic bias due to the penalty for each nonzero parameter is $a_{n1}$. Since we penalized only on the off-diagonal elements, the total bias induced by the penalty is asymptotically of order $s_{n1}a_{n1}$. The square of this total bias over all nonzero elements is of order $O_P\{(p_n + s_{n1})\log p_n/n\}$ under condition (B).

Theorem 1 states explicitly how the number of nonzero elements and dimensionality affect the rate of convergence. Since there are $(p_n + s_{n1})$ nonzero elements and each of them can be estimated at best with rate $n^{-1/2}$, the total square errors are at least of rate $(p_n + s_{n1})/n$. The price that we pay for high-dimensionality is merely a logarithmic factor $\log p_n$. The results holds as long as $(p_n + s_{n1})/n$ is at a rate $O((\log p_n)^{-k})$ with some $k > 1$, which decays to zero slowly. This means that in practice $p_n$ can be comparable to $n$ without violating the results. The condition here is not minimum possible; we expect it holds for $p \gg n$. Here, we refer the local minimizer as an interior point within a given close set such that it minimizes the target function. Following a similar argument to Huang, Horowitz and Ma (2008), the local minimizer in Theorem 1 can be taken as the global minimizer with additional conditions on the tail of the penalty function.

Theorem 1 is also applicable to the $L_1$-penalty function, where the local minimizer for $\lambda_{n1}$ can be relaxed. In this case, the local minimizer becomes the global minimizer. The asymptotic bias of the $L_1$-penalized estimate is given in the term $s_{n1}a_{n1} = s_{n1}\lambda_{n1}$ as shown in the technical proof. In order to control the bias, we impose condition (B), which entails an upper bound on $\lambda_{n1} = O(\{(1 + p_n/(s_{n1} + 1))\log p_n/n\}^{1/2})$. The bias problem due to the $L_1$-penalty for finite parameter has already been unveiled by Fan and Li (2001) and Zou (2006).

Next, we show the sparsistency of the penalized estimator from (1.1). We use $S^c$ to denote the complement of a set $S$.

THEOREM 2 (Sparsistency). *Under the conditions given in Theorem 1, for any local minimizer of (1.1) satisfying $\|\hat{\mathbf{\Omega}} - \mathbf{\Omega}_0\|_F^2 = O_P\{(p_n + s_{n1})\log p_n/n\}$ and $\|\hat{\mathbf{\Omega}} - \mathbf{\Omega}_0\|^2 = O_P(\eta_n)$ for a sequence of $\eta_n \to 0$, if $\log p_n/n + \eta_n = O(\lambda_{n1}^2)$, then with probability tending to 1, $\hat{\omega}_{ij} = 0$ for all $(i,j) \in S_1^c$.*

First, since $\|M\|^2 \leq \|M\|_F^2$ for any matrix $M$, we can always take $\eta_n = (p_n + s_{n1})\log p_n/n$ in Theorem 2, but this will result in more stringent requirement on the number of zero elements when $L_1$-penalty is used, as we now explain. The sparsistency requires a lower bound on the rate of the regularization parameter $\lambda_{n1}$. On the other hand, condition (B) imposes an upper bound on $\lambda_{n1}$ when $L_1$-penalty is used in order to control the biases. Explicitly, we need, for $L_1$-penalized likelihood,

$$(2.3) \qquad \log p_n/n + \eta_n = O(\lambda_{n1}^2) = (1 + p_n/(s_{n1}+1))\log p_n/n$$

for both consistency and sparsistency to be satisfied. We present two scenarios here for the two bounds to be compatible, making use of the inequalities $\|M\|_F^2/p_n \leq \|M\|^2 \leq \|M\|_F^2$ for a matrix $M$ of size $p_n$.

1. We always have $\|\hat{\boldsymbol{\Omega}} - \boldsymbol{\Omega}_0\| \leq \|\hat{\boldsymbol{\Omega}} - \boldsymbol{\Omega}_0\|_F$. In the worst case scenario where they have the same order, $\|\hat{\boldsymbol{\Omega}} - \boldsymbol{\Omega}_0\|^2 = O_P((p_n + s_{n1})\log p_n/n)$, so that $\eta_n = (p_n + s_{n1})\log p_n/n$. It is then easy to see from (2.3) that the two bounds are compatible only when $s_{n1} = O(1)$.
2. We also have $\|\hat{\boldsymbol{\Omega}} - \boldsymbol{\Omega}_0\|_F^2/p_n \leq \|\hat{\boldsymbol{\Omega}} - \boldsymbol{\Omega}_0\|^2$. In the optimistic scenario where they have the same order,

$$\|\hat{\boldsymbol{\Omega}} - \boldsymbol{\Omega}_0\|^2 = O_P((1 + s_{n1}/p_n)\log p_n/n).$$

Hence, $\eta_n = (1 + s_{n1}/p_n)\log p_n/n$, and compatibility of the bounds requires $s_{n1} = O(p_n)$.

Hence, even in the optimistic scenario, consistency and sparsistency are guaranteed only when $s_{n1} = O(p_n)$ if the $L_1$-penalty is used, that is, the precision matrix has to be sparse enough.

However, if the penalty function used is unbiased, like the SCAD or the hard-thresholding penalty, we do not impose an extra upper bound for $\lambda_{n1}$ since its first derivative $p'_{\lambda_{n1}}(|\theta|)$ goes to zero fast enough as $|\theta|$ increases [exactly equals zero for the SCAD and hard-thresholding penalty functions, when $n$ is sufficiently large; see condition (B) and the explanation thereof]. Thus, $\lambda_{n1}$ is allowed to decay to zero slowly, allowing even the largest order $s_{n1} = O(p_n^2)$.

We remark that asymptotic normality for the estimators of the elements in $S_1$ have been established in a previous version of this paper. We omit it here for brevity.

2.5. *Properties of sparse inverse correlation matrix estimation.* The inverse correlation matrix $\boldsymbol{\Psi}_0$ retains the same sparsity structure of $\boldsymbol{\Omega}_0$. Consistency and sparsistency results can be achieved with $p_n$ as large as $\log p_n = o(n)$, as long as $(s_{n1}+1)(\log p_n)^k/n = O(1)$ for some $k > 1$ as $n \to \infty$. We minimize, w.r.t. $\boldsymbol{\Psi} = (\psi_{ij})$,

$$(2.4) \qquad \operatorname{tr}(\boldsymbol{\Psi}\hat{\boldsymbol{\Gamma}}_{\mathbf{S}}) - \log|\boldsymbol{\Psi}| + \sum_{i \neq j} p_{\nu_{n1}}(|\psi_{ij}|),$$

where $\hat{\boldsymbol{\Gamma}}_{\mathbf{S}} = \hat{\mathbf{W}}^{-1}\mathbf{S}\hat{\mathbf{W}}^{-1}$ is the sample correlation matrix, with $\hat{\mathbf{W}}^2 = \mathbf{D}_{\mathbf{S}}$ being the diagonal matrix with diagonal elements of $\mathbf{S}$, and $\nu_{n1}$ is a regularization parameter. After obtaining $\hat{\boldsymbol{\Psi}}$, $\boldsymbol{\Omega}_0$ can also be estimated by $\tilde{\boldsymbol{\Omega}} = \hat{\mathbf{W}}^{-1}\hat{\boldsymbol{\Psi}}\hat{\mathbf{W}}^{-1}$.

To present the rates of convergence for $\hat{\boldsymbol{\Psi}}$ and $\tilde{\boldsymbol{\Omega}}$, we define

$$c_{n1} = \max_{(i,j)\in S_1} p'_{\nu_{n1}}(|\psi^0_{ij}|), \qquad d_{n1} = \max_{(i,j)\in S_1} p''_{\nu_{n1}}(|\psi^0_{ij}|),$$

where $\boldsymbol{\Psi}_0 = (\psi^0_{ij})$ and modify condition (D) to (D') with $\lambda_{n1}$ there replaced by $\nu_{n1}$, and impose

(B') $c_{n1} = O(\{\log p_n/n\}^{1/2})$, $d_{n1} = o(1)$. Also, $\min_{(i,j)\in S_1} |\psi^0_{ij}|/\nu_{n1} \to \infty$ as $n \to \infty$.

THEOREM 3. *Under regularity conditions* (A), (B'), (C) *and* (D'), *if* $(s_{n1}+1)(\log p_n)^k/n = O(1)$ *for some* $k > 1$ *and* $(s_{n1}+1)\log p_n/n = o(\nu_{n1}^2)$, *then there exists a local minimizer* $\hat{\boldsymbol{\Psi}}$ *for (2.4) such that* $\|\hat{\boldsymbol{\Psi}} - \boldsymbol{\Psi}_0\|_F^2 = O_P(s_{n1}\log p_n/n)$ *and* $\|\tilde{\boldsymbol{\Omega}} - \boldsymbol{\Omega}_0\|^2 = O_P((s_{n1}+1)\log p_n/n)$ *under the operator norm. For the $L_1$-penalty, we only need* $\log p_n/n = O(\nu_{n1}^2)$.

Note that we can allow $p_n \gg n$ without violating the result as long as $\log p_n/n = o(1)$. Note also that an order of $\{p_n \log p_n/n\}^{1/2}$ is removed by estimating the inverse correlation rather than the precision matrix, which is somewhat surprising since the inverse correlation matrix, unlike the correlation matrix, does not have known diagonal elements that contribute no errors to the estimation. This can be explained and proved as follows. If $s_{n1} = O(p_n)$, the result is obvious. When $s_{n1} = o(p_n)$, most of the off-diagonal elements are zero. Indeed, there are at most $O(s_{n1})$ columns of the inverse correlation matrix which contain at least one nonzero element. The rest of the columns that have all zero off-diagonal elements must have diagonal entries 1. These columns represent variables that are actually uncorrelated from the rest. Now, it is easy to see from (2.4) that these diagonal elements, which are one, are all estimated exactly as one with no estimation error. Hence, an order of $(p_n \log p_n/n)^{1/2}$ is not present even in the case of estimating the inverse correlation matrix.

For the $L_1$-penalty, our result reduces to that given in Rothman et al. (2008). We offer the sparsistency result as follows.

THEOREM 4 (Sparsistency). *Under the conditions given in Theorem 3, for any local minimizer of (2.4) satisfying* $\|\hat{\boldsymbol{\Psi}} - \boldsymbol{\Psi}_0\|_F^2 = O_P(s_{n1}\log p_n/n)$ *and* $\|\hat{\boldsymbol{\Psi}} - \boldsymbol{\Psi}_0\|^2 = O_P(\eta_n)$ *for some* $\eta_n \to 0$, *if* $\log p_n/n + \eta_n = O(\nu_{n1}^2)$, *then with probability tending to 1*, $\hat{\psi}_{ij} = 0$ *for all* $(i,j) \in S_1^c$.

The proof follows exactly the same as that for Theorem 2 in Section 2.4, and is thus omitted.

For the $L_1$-penalty, control of bias and sparsistency require $\nu_{n1}$ to satisfy bounds like (2.3):

(2.5) $$\log p_n/n + \eta_n = O(\nu_{n1}^2) = \log p_n/n.$$

This leads to two scenarios:

1. The worst case scenario has

$$\|\hat{\boldsymbol{\Psi}} - \boldsymbol{\Psi}_0\|^2 = \|\hat{\boldsymbol{\Psi}} - \boldsymbol{\Psi}_0\|_F^2 = O_P(s_{n1} \log p_n/n),$$

meaning $\eta_n = s_{n1} \log p_n/n$. Then compatibility of the bounds in (2.5) requires $s_{n1} = O(1)$.

2. The optimistic scenario has

$$\|\hat{\boldsymbol{\Psi}} - \boldsymbol{\Psi}_0\|^2 = \|\hat{\boldsymbol{\Psi}} - \boldsymbol{\Psi}_0\|_F^2/p_n = O_P(s_{n1}/p_n \cdot \log p_n/n),$$

meaning $\eta_n = s_{n1}/p_n \cdot \log p_n/n$. Then compatibility of the bounds in (2.5) requires $s_{n1} = O(p_n)$.

On the other hand, for penalties like the SCAD or the hard-thresholding penalty, we do not need an upper bound for $s_{n1}$. Hence, we only need $(s_{n1} + 1)(\log p_n)^k/n = O(1)$ as $n \to \infty$ for some $k > 1$. It is clear that SCAD results in better sampling properties than the $L_1$-penalized estimator in precision or inverse correlation matrix estimation.

**3. Estimation of sparse covariance matrix.** In this section, we analyze the sparse covariance matrix estimation using the penalized likelihood (1.3). Then it is modified to estimating the correlation matrix, which improves the rate of convergence. We assume that the $\mathbf{y}_i$'s are i.i.d. $N(\mathbf{0}, \boldsymbol{\Sigma}_0)$ throughout this section.

3.1. *Properties of sparse covariance matrix estimation.* Let $S_2 = \{(i,j): \sigma_{ij}^0 \neq 0\}$, where $\boldsymbol{\Sigma}_0 = (\sigma_{ij}^0)$. Denote $s_{n2} = |S_2| - p_n$, so that $s_{n2}$ is the number of nonzero elements in $\boldsymbol{\Sigma}_0$ on the off-diagonal entries. Put

$$a_{n2} = \max_{(i,j) \in S_2} p'_{\lambda_{n2}}(|\sigma_{ij}^0|), \qquad b_{n2} = \max_{(i,j) \in S_2} p''_{\lambda_{n2}}(|\sigma_{ij}^0|).$$

Technical conditions in Section 2 need some revision. In particular, condition (D) now becomes condition (D2) with $\lambda_{n1}$ there replaced by $\lambda_{n2}$. Condition (B) should now be

(B2) $a_{n2} = O(\{(1+p_n/(s_{n2}+1))\log p_n/n\}^{1/2})$, $b_{n2} = o(1)$, and $\min_{(i,j) \in S_2} |\sigma_{ij}^0|/\lambda_{n2} \to \infty$ as $n \to \infty$.

THEOREM 5 (Rate of convergence). *Under regularity conditions* (A), (B2), (C) *and* (D2), *if* $(p_n + s_{n2})(\log p_n)^k/n = O(1)$ *for some* $k > 1$ *and* $(p_n + s_{n2}) \log p_n/n = O(\lambda_{n2}^2)$, *then there exists a local minimizer* $\hat{\boldsymbol{\Sigma}}$ *such that* $\|\hat{\boldsymbol{\Sigma}} - \boldsymbol{\Sigma}_0\|_F^2 = O_P((p_n + s_{n2}) \log p_n/n)$. *For the $L_1$-penalty, we only need* $\log p_n/n = O(\lambda_{n2}^2)$.

When the $L_1$-penalty is used, the condition for $\lambda_{n2}$ is relaxed to $\log p_n/n = O(\lambda_{n2}^2)$. Like the case for precision matrix estimation, the asymptotic bias due to the $L_1$-penalty is of order $s_{n2}a_{n2} = s_{n2}\lambda_{n2}$. To control this term, for the $L_1$-penalty, we require $\lambda_{n2} = O(\{(1 + p_n/(s_{n2}+1))\log p_n/n\}^{1/2})$.

THEOREM 6 (Sparsistency). *Under the conditions given in Theorem 5, for any local minimizer $\hat{\boldsymbol{\Sigma}}$ of (1.3) satisfying $\|\hat{\boldsymbol{\Sigma}} - \boldsymbol{\Sigma}_0\|_F^2 = O_P((p_n + s_{n2})\log p_n/n)$ and $\|\hat{\boldsymbol{\Sigma}} - \boldsymbol{\Sigma}_0\|^2 = O_P(\eta_n)$ for some $\eta_n \to 0$, if $\log p_n/n + \eta_n = O(\lambda_{n2}^2)$, then with probability tending to 1, $\hat{\sigma}_{ij} = 0$ for all $(i,j) \in S_2^c$.*

For the $L_1$-penalized likelihood, controlling of bias for consistency together with sparsistency requires

$$(3.1) \qquad \log p_n/n + \eta_n = O(\lambda_{n2}^2) = (1 + p_n/(s_{n2}+1))\log p_n/n.$$

This is the same condition as (2.3), and hence in the worst case scenario where

$$\|\hat{\boldsymbol{\Sigma}} - \boldsymbol{\Sigma}_0\|^2 = \|\hat{\boldsymbol{\Sigma}} - \boldsymbol{\Sigma}_0\|_F^2 = O_P((p_n + s_{n2})\log p_n/n),$$

we need $s_{n2} = O(1)$. In the optimistic scenario where

$$\|\hat{\boldsymbol{\Sigma}} - \boldsymbol{\Sigma}_0\|^2 = \|\hat{\boldsymbol{\Sigma}} - \boldsymbol{\Sigma}_0\|_F^2/p_n,$$

we need $s_{n2} = O(p_n)$. In both cases, the matrix $\boldsymbol{\Sigma}_0$ has to be very sparse, but the former is much sparser.

On the other hand, if unbiased penalty functions like the SCAD or hard-thresholding penalty are used, we do not need an upper bound on $\lambda_{n2}$ since the bias $a_{n2} = 0$ for sufficiently large $n$. This gives more flexibility on the order of $s_{n2}$.

Similar to Section 2, asymptotic normality for the estimators of the elements in $S_2$ can be proved under certain assumptions.

3.2. *Properties of sparse correlation matrix estimation.* The correlation matrix $\boldsymbol{\Gamma}_0$ retains the same sparsity structure of $\boldsymbol{\Sigma}_0$ with known diagonal elements. This special structure allows us to estimate $\boldsymbol{\Gamma}_0$ more accurately. To take advantage of the known diagonal elements, the sparse correlation matrix $\boldsymbol{\Gamma}_0$ is estimated by minimizing w.r.t. $\boldsymbol{\Gamma} = (\gamma_{ij})$,

$$(3.2) \qquad \mathrm{tr}(\boldsymbol{\Gamma}^{-1}\hat{\boldsymbol{\Gamma}}_{\mathbf{S}}) + \log|\boldsymbol{\Gamma}| + \sum_{i \neq j} p_{\nu_{n2}}(|\gamma_{ij}|),$$

where $\nu_{n2}$ is a regularization parameter. After obtaining $\hat{\boldsymbol{\Gamma}}$, $\boldsymbol{\Sigma}_0$ can be estimated by $\tilde{\boldsymbol{\Sigma}} = \hat{\mathbf{W}}\hat{\boldsymbol{\Gamma}}\hat{\mathbf{W}}$.

To present the rates of convergence for $\hat{\boldsymbol{\Gamma}}$ and $\tilde{\boldsymbol{\Sigma}}$, we define

$$c_{n2} = \max_{(i,j) \in S_2} p'_{\nu_{n2}}(|\gamma_{ij}^0|), \qquad d_{n2} = \max_{(i,j) \in S_2} p''_{\nu_{n2}}(|\gamma_{ij}^0|),$$

where $\mathbf{\Gamma}_0 = (\gamma_{ij}^0)$. We modify condition (D) to (D2′) with $\lambda_{n2}$ there replaced by $\nu_{n2}$, and (B) to (B′) as follows:

(B2′) $c_{n2} = O(\{\log p_n/n\}^{1/2})$, $d_{n2} = o(1)$, and $\min_{(i,j) \in S_2} |\gamma_{ij}^0|/\nu_{n2} \to \infty$ as $n \to \infty$.

THEOREM 7. *Under regularity conditions (A), (B2′), (C) and (D2′), if $(p_n + s_{n2})(\log p_n)^k/n = O(1)$ for some $k > 1$ and $(s_{n2} + 1)\log p_n/n = O(\nu_{n2}^2)$, then there exists a local minimizer $\hat{\mathbf{\Gamma}}$ for (3.2) such that*

$$\|\hat{\mathbf{\Gamma}} - \mathbf{\Gamma}_0\|_F^2 = O_P(s_{n2} \log p_n/n).$$

*In addition, for the operator norm, we have*

$$\|\tilde{\mathbf{\Sigma}} - \mathbf{\Sigma}_0\|^2 = O_P\{(s_{n2} + 1)\log p_n/n\}.$$

*For the $L_1$-penalty, we only need $\log p_n/n = O(\nu_{n2}^2)$.*

The proof is sketched in Section 5. This theorem shows that the correlation matrix, like the inverse correlation matrix, can be estimated more accurately, since diagonal elements are known to be one.

THEOREM 8 (Sparsistency). *Under the conditions given in Theorem 7, for any local minimizer $\hat{\mathbf{\Gamma}}$ of (3.2) satisfying $\|\hat{\mathbf{\Gamma}} - \mathbf{\Gamma}_0\|_F^2 = O_P(s_{n2} \log p_n/n)$ and $\|\hat{\mathbf{\Gamma}} - \mathbf{\Gamma}_0\|^2 = O_P(\eta_n)$ for some $\eta_n \to 0$, if $\log p_n/n + \eta_n = O(\nu_{n2}^2)$, then with probability tending to 1, $\hat{\gamma}_{ij} = 0$ for all $(i,j) \in S_2^c$.*

The proof follows exactly the same as that of Theorem 6 in Section 5, and is omitted. For the $L_1$-penalized likelihood, controlling of bias and sparsistency requires

(3.3) $$\log p_n/n + \eta_n = O(\nu_{n2}^2) = \log p_n/n.$$

This is the same condition as (2.5), hence in the worst scenario where

$$\|\hat{\mathbf{\Gamma}} - \mathbf{\Gamma}_0\|^2 = \|\hat{\mathbf{\Gamma}} - \mathbf{\Gamma}_0\|_F^2 = O_P(s_{n2} \log p_n/n),$$

we need $s_{n2} = O(1)$. In the optimistic scenario where

$$\|\hat{\mathbf{\Gamma}} - \mathbf{\Gamma}_0\|^2 = \|\hat{\mathbf{\Gamma}} - \mathbf{\Gamma}_0\|_F^2/p_n = O_P(s_{n2}/p_n \cdot \log p_n/n),$$

we need $s_{n2} = O(p_n)$.

The use of unbiased penalty functions like the SCAD or the hard-thresholding penalty, similar to results in the previous sections, does not impose an upper bound on the regularization parameter since bias $c_{n2} = 0$ for sufficiently large $n$. This gives more flexibility to the order of $s_{n2}$ allowed.

**4. Extension to sparse Cholesky decomposition.** Pourahmadi ([1999](#)) proposed the modified Cholesky decomposition (MCD) which facilitates the sparse estimation of $\boldsymbol{\Omega}$ through penalization. The idea is to represent zero-mean data $\mathbf{y} = (y_1, \ldots, y_{p_n})^T$ using the autoregressive model:

$$(4.1) \qquad y_i = \sum_{j=1}^{i-1} \phi_{ij} y_j + \epsilon_i \quad \text{and} \quad \mathbf{T}\boldsymbol{\Sigma}\mathbf{T}^T = \mathbf{D},$$

where $\mathbf{T}$ is the unique unit lower triangular matrix with ones on its diagonal and $(i,j)$th element being $-\phi_{ij}$ for $j < i$, and $\mathbf{D}$ is diagonal with $i$th element being $\sigma_i^2 = \text{var}(\epsilon_i)$. The optimization problem is unconstrained (since the $\phi_{ij}$'s are free variables), and the estimate for $\boldsymbol{\Omega}$ is always positive-definite.

Huang et al. ([2006](#)) and Levina, Rothman and Zhu ([2008](#)) both used the MCD for estimating $\boldsymbol{\Omega}_0$. The former maximized the log-likelihood (ML) over $\mathbf{T}$ and $\mathbf{D}$ simultaneously, while the latter suggested also a least square version (LS), with $\mathbf{D}$ being first set to the identity matrix and then minimizing over $\mathbf{T}$ to obtain $\hat{\mathbf{T}}$. The latter corresponds to the original Cholesky decomposition. The sparse Cholesky factor can be estimated through minimizing

$$(4.2) \quad \text{(ML): } q_3(\mathbf{T}, \mathbf{D}) = \text{tr}(\mathbf{T}^T \mathbf{D}^{-1} \mathbf{T} \mathbf{S}) + \log|\mathbf{D}| + 2\sum_{i<j} p_{\lambda_{n3}}(|t_{ij}|).$$

This is indeed the same as (1.1) with the substitution of $\boldsymbol{\Omega} = \mathbf{T}^T \mathbf{D}^{-1} \mathbf{T}$ and penalization parameter $\lambda_{n3}$. Noticing that (4.1) can be written as $\mathbf{T}\mathbf{y} = \boldsymbol{\varepsilon}$, the least square version is to minimize $\text{tr}(\boldsymbol{\varepsilon}\boldsymbol{\varepsilon}^T) = \text{tr}(\mathbf{T}^T \mathbf{T} \mathbf{y}\mathbf{y}^T)$ in the matrix notation. Aggregating the $n$ observations and adding penalty functions, the least-square criterion is to minimize

$$(4.3) \qquad \text{(LS): } q_4(\mathbf{T}) = \text{tr}(\mathbf{T}^T \mathbf{T} \mathbf{S}) + 2\sum_{i<j} p_{\lambda_{n4}}(|t_{ij}|).$$

In view of the results in Sections 2.5 and 3.2, we can also write the sample covariance matrix in (4.2) as $\mathbf{S} = \hat{\mathbf{W}}\hat{\boldsymbol{\Gamma}}_{\mathbf{S}}\hat{\mathbf{W}}$ and then replace $\mathbf{D}^{-1/2}\mathbf{T}\hat{\mathbf{W}}$ by $\mathbf{T}$, resulting in the normalized (NL) version as follows:

$$(4.4) \qquad \text{(NL): } q_5(\mathbf{T}) = \text{tr}(\mathbf{T}^T \mathbf{T} \hat{\boldsymbol{\Gamma}}_{\mathbf{S}}) - 2\log|\mathbf{T}| + 2\sum_{i<j} p_{\lambda_{n5}}(|t_{ij}|).$$

We will also assume the $\mathbf{y}_i$'s are i.i.d. $N(\mathbf{0}, \boldsymbol{\Sigma}_0)$ as in the last section.

4.1. *Properties of sparse Cholesky factor estimation.* Since all the $\mathbf{T}$'s introduced in the three models above have the same sparsity structure, let $S$ and $s_{n3}$ be the nonzero set and number of nonzeros associated with each $\mathbf{T}$ above. Define

$$a_{n3} = \max_{(i,j)\in S} p'_{\lambda_{n3}}(|t^0_{ij}|), \qquad b_{n3} = \max_{(i,j)\in S} p''_{\lambda_{n3}}(|t^0_{ij}|).$$

For (ML), condition (D) is adapted to (D3) with $\lambda_{n1}$ there replaced by $\lambda_{n3}$. Condition (B) is modified as:

(B3) $a_{n3} = O(\{(1+p_n/(s_{n3}+1))\log p_n/n\}^{1/2})$, $b_{n3} = o(1)$ and $\min_{(i,j)\in S} |\phi_{ij}^0|/\lambda_{n3} \to \infty$ as $n \to \infty$.

After obtaining $\hat{\mathbf{T}}$ and $\hat{\mathbf{D}}$ from minimizing (ML), we set $\hat{\mathbf{\Omega}} = \hat{\mathbf{T}}^T \hat{\mathbf{D}}^{-1} \hat{\mathbf{T}}$.

THEOREM 9. *Under regularity conditions* (A), (B3), (C), (D3), *if* $(p_n + s_{n3})(\log p_n)^k/n = O(1)$ *for some* $k > 1$ *and* $(p_n + s_{n3})\log p_n/n = O(\lambda_{n3}^2)$, *then there exists a local minimizer* $\hat{\mathbf{T}}$ *and* $\hat{\mathbf{D}}$ *for* (ML) *such that* $\|\hat{\mathbf{T}} - \mathbf{T}_0\|_F^2 = O_P(s_{n3}\log p_n/n)$, $\|\hat{\mathbf{D}} - \mathbf{D}_0\|_F^2 = O_P(p_n\log p_n/n)$ *and* $\|\hat{\mathbf{\Omega}} - \mathbf{\Omega}_0\|_F^2 = O_P\{(p_n+s_{n3})\log p_n/n\}$. *For the* $L_1$-*penalty, we only need* $\log p_n/n = O(\lambda_{n3}^2)$.

The proof is similar to those of Theorems 5 and 7 and is omitted. The Cholesky factor $\mathbf{T}$ has ones on its main diagonal without the need for estimation. Hence, the rate of convergence is faster than $\hat{\mathbf{\Omega}}$.

THEOREM 10 (Sparsistency). *Under the conditions in Theorem 9, for any local minimizer* $\hat{\mathbf{T}}$, $\hat{\mathbf{D}}$ *of* (4.2) *satisfying* $\|\hat{\mathbf{T}} - \mathbf{T}_0\|_F^2 = O_P(s_{n3}\log p_n/n)$ *and* $\|\hat{\mathbf{D}} - \mathbf{D}_0\|_F^2 = O_P(p_n\log p_n/n)$, *if* $\log p_n/n + \eta_n + \zeta_n = O(\lambda_{n3}^2)$, *then sparsistency holds for* $\hat{\mathbf{T}}$, *provided that* $\|\hat{\mathbf{T}} - \mathbf{T}_0\|^2 = O_P(\eta_n)$ *and* $\|\hat{\mathbf{D}} - \mathbf{D}_0\|^2 = O_P(\zeta_n)$, *for some* $\eta_n, \zeta_n \to 0$.

The proof is in Section 5. For the $L_1$-penalized likelihood, control of bias and sparsistency impose the following:

(4.5) $\quad \log p_n/n + \eta_n + \zeta_n = O(\lambda_{n3}^2) = (1 + p_n/(s_{n3}+1))\log p_n/n$.

The worst scenario corresponds to $\eta_n = s_{n3}\log p_n/n$ and $\zeta_n = p_n\log p_n/n$, so that we need $s_{n3} = O(1)$. The optimistic scenario corresponds to $\eta_n = s_{n3}/p_n \cdot \log p_n/n$ and $\zeta_n = \log p_n/n$, so that we need $s_{n3} = O(p_n)$.

On the other hand, such a restriction is not needed for unbiased penalties like the SCAD or hard-thresholding penalty, giving more flexibility on the order of $s_{n3}$.

4.2. *Properties of sparse normalized Cholesky factor estimation.* We now turn to analyzing the normalized penalized likelihood (4.4). With $\mathbf{T} = (t_{ij})$ in (NL) which is lower triangular, define

$$a_{n5} = \max_{(i,j)\in S} p'_{\lambda_{n5}}(|t_{ij}^0|), \qquad b_{n5} = \max_{(i,j)\in S} p''_{\lambda_{n5}}(|t_{ij}^0|).$$

Condition (D) is now changed to (D5) with $\lambda_{n1}$ there replaced by $\lambda_{n5}$. Condition (B) is now substituted by

(B5) $a_{n5}^2 = O(\log p_n/n)$, $b_{n5} = o(1)$, $\min_{(i,j)\in S} |t_{ij}^0|/\lambda_{n5} \to \infty$ as $n \to \infty$.

THEOREM 11 (Rate of convergence). *Under regularity conditions* (A), (B5), (C) *and* (D5), *if* $s_{n3}(\log p_n)^k/n = O(1)$ *for some* $k > 1$ *and* $(s_{n3} + 1)\log p_n/n = O(\lambda_{n5}^2)$, *then there exists a local minimizer* $\hat{\mathbf{T}}$ *for* (NL) *such that* $\|\hat{\mathbf{T}} - \mathbf{T}_0\|_F^2 = O_P(s_{n3}\log p_n/n)$ *and rate of convergence in the Frobenius norm*

$$\|\hat{\mathbf{\Omega}} - \mathbf{\Omega}_0\|_F^2 = O_P\{(p_n + s_{n3})\log p_n/n\}$$

*and in the operator norm, it is improved to*

$$\|\hat{\mathbf{\Omega}} - \mathbf{\Omega}_0\|^2 = O_P\{(s_{n3} + 1)\log p_n/n\}.$$

*For the $L_1$-penalty, we only need* $\log p_n/n = O(\lambda_{n5}^2)$.

The proof is similar to that of Theorems 5 and 7 and is omitted. In this theorem, like Lemma 3, we can have $p_n$ so that $p_n/n$ goes to a constant less than 1. It is evident that normalizing with $\hat{\mathbf{W}}$ results in an improvement in the rate of convergence in operator norm.

THEOREM 12 (Sparsistency). *Under the conditions given in Theorem 11, for any local minimizer* $\hat{\mathbf{T}}$ *of (4.4) satisfying* $\|\hat{\mathbf{T}} - \mathbf{T}_0\|_F^2 = O_P(s_{n3}\log p_n/n)$ *if* $\log p_n/n + \eta_n = O(\lambda_{n5}^2)$, *then sparsistency holds for* $\hat{\mathbf{T}}$, *provided that* $\|\hat{\mathbf{T}} - \mathbf{T}_0\|^2 = O(\eta_n)$ *for some* $\eta_n \to 0$.

Proof is omitted since it goes exactly the same as that of Theorem 10. The above results apply also to the $L_1$-penalized estimator. For simultaneous persistency and optimal rate of convergence using the $L_1$-penalty, the biases inherent in it induce the restriction $s_{n3} = O(1)$ in the worst scenario where $\eta_n^2 = s_{n3}\log p_n/n$, and $s_{n3} = O(p_n)$ in the optimistic scenario where $\eta_n^2 = s_{n3}/p_n \cdot \log p_n/n$. This restriction does not apply to the SCAD and other asymptotically unbiased penalty functions.

**5. Proofs.** We first prove three lemmas. The first one concerns with inequalities involving the operator and the Frobenius norms. The other two concern with order estimation for elements in a matrix of the form $\mathbf{A}(\mathbf{S} - \mathbf{\Sigma}_0)\mathbf{B}$, which are useful in proving results concerning sparsistency.

LEMMA 1. *Let* $\mathbf{A}$ *and* $\mathbf{B}$ *be real matrices such that the product* $\mathbf{AB}$ *is defined. Then, defining* $\|\mathbf{A}\|_{\min}^2 = \lambda_{\min}(\mathbf{A}^T\mathbf{A})$, *we have*

(5.1) $$\|\mathbf{A}\|_{\min}\|\mathbf{B}\|_F \leq \|\mathbf{AB}\|_F \leq \|\mathbf{A}\|\|\mathbf{B}\|_F.$$

*In particular, if* $\mathbf{A} = (a_{ij})$, *then* $|a_{ij}| \leq \|\mathbf{A}\|$ *for each* $i, j$.

PROOF. Write $\mathbf{B} = (\mathbf{b}_1, \ldots, \mathbf{b}_q)$, where $\mathbf{b}_i$ is the $i$th column vector in $\mathbf{B}$. Then

$$\|\mathbf{AB}\|_F^2 = \operatorname{tr}(\mathbf{B}^T\mathbf{A}^T\mathbf{AB}) = \sum_{i=1}^q \mathbf{b}_i^T\mathbf{A}^T\mathbf{A}\mathbf{b}_i \leq \lambda_{\max}(\mathbf{A}^T\mathbf{A})\sum_{i=1}^q \|\mathbf{b}_i\|^2$$
$$= \|\mathbf{A}\|^2 \|\mathbf{B}\|_F^2.$$

Similarly,

$$\|\mathbf{AB}\|_F^2 = \sum_{i=1}^q \mathbf{b}_i^T\mathbf{A}^T\mathbf{A}\mathbf{b}_i \geq \lambda_{\min}(\mathbf{A}^T\mathbf{A})\sum_{i=1}^q \|\mathbf{b}_i\|^2$$
$$= \|\mathbf{A}\|_{\min}^2 \|\mathbf{B}\|_F^2,$$

which completes the proof of (5.1). To prove $|a_{ij}| \leq \|\mathbf{A}\|$, note that $a_{ij} = \mathbf{e}_i^T\mathbf{A}\mathbf{e}_j$, where $\mathbf{e}_i$ is the unit column vector with one at the $i$th position, and zero elsewhere. Hence, using (5.1),

$$|a_{ij}| = |\mathbf{e}_i^T\mathbf{A}\mathbf{e}_j| \leq \|\mathbf{A}\mathbf{e}_j\|_F \leq \|\mathbf{A}\| \cdot \|\mathbf{e}_j\|_F = \|\mathbf{A}\|,$$

and this completes the proof of the lemma. □

LEMMA 2. *Let $\mathbf{S}$ be a sample covariance matrix of a random sample $\{\mathbf{y}_i\}_{1 \leq i \leq n}$, with $E(\mathbf{y}_i) = 0$ and $\operatorname{var}(\mathbf{y}_i) = \mathbf{\Sigma}_0$. Let $\mathbf{y}_i = (y_{i1}, \ldots, y_{ip_n})$ with $y_{ij} \sim F_j$, where $F_j$ is the c.d.f. of $y_{ij}$, and let $G_j$ be the c.d.f. of $y_{ij}^2$, with*

(5.2) $$\max_{1 \leq i \leq p_n} \int_0^\infty \exp(\lambda t)\, dG_j(t) < \infty, \qquad 0 < |\lambda| < \lambda_0,$$

*for some $\lambda_0 > 0$. Assume $\log p_n/n = o(1)$, and that $\mathbf{\Sigma}_0$ has eigenvalues uniformly bounded above as $n \to \infty$. Then for constant matrices $\mathbf{A}$ and $\mathbf{B}$ with $\|\mathbf{A}\|, \|\mathbf{B}\| = O(1)$, we have $\max_{i,j} |(\mathbf{A}(\mathbf{S} - \mathbf{\Sigma}_0)\mathbf{B})_{ij}| = O_P(\{\log p_n/n\}^{1/2})$.*

REMARK. The conditions on the $y_{ij}$'s above are the same as those used in Bickel and Levina (2008b) for relaxing the normality assumption.

PROOF OF LEMMA 2. Let $\mathbf{x}_i = \mathbf{A}\mathbf{y}_i$ and $\mathbf{w}_i = \mathbf{B}^T\mathbf{y}_i$. Define $\mathbf{u}_i = (\mathbf{x}_i^T, \mathbf{w}_i^T)^T$, with covariance matrix

$$\mathbf{\Sigma}_{\mathbf{u}} = \operatorname{var}(\mathbf{u}_i) = \begin{pmatrix} \mathbf{A}\mathbf{\Sigma}_0\mathbf{A}^T & \mathbf{A}\mathbf{\Sigma}_0\mathbf{B} \\ \mathbf{B}^T\mathbf{\Sigma}_0\mathbf{A}^T & \mathbf{B}^T\mathbf{\Sigma}_0\mathbf{B} \end{pmatrix}.$$

Since $\|(\mathbf{A}^T\mathbf{B})^T\| \leq (\|\mathbf{A}\|^2 + \|\mathbf{B}\|^2)^{1/2} = O(1)$ and $\|\mathbf{\Sigma}_0\| = O(1)$ uniformly, we have $\|\mathbf{\Sigma}_{\mathbf{u}}\| = O(1)$ uniformly. Then, with $\mathbf{S}_{\mathbf{u}} = n^{-1}\sum_{i=1}^n \mathbf{u}_i\mathbf{u}_i^T$, which is the sample covariance matrix for the random sample $\{\mathbf{u}_i\}_{1 \leq i \leq n}$, by Lemma

A.3 of Bickel and Levina (2008b) which holds under the assumption for the $y_{ij}$'s and $\log p_n/n = o(1)$, we have

$$\max_{i,j}|(\mathbf{S_u} - \boldsymbol{\Sigma_u})_{ij}| = O_P(\{\log p_n/n\}^{1/2}).$$

In particular, it means that

$$\max_{i,j}|(\mathbf{A}(\mathbf{S} - \boldsymbol{\Sigma}_0)\mathbf{B})_{ij}| = \left(n^{-1}\sum_{r=1}^{n}\mathbf{x}_r\mathbf{w}_r^T - \mathbf{A}\boldsymbol{\Sigma}_0\mathbf{B}\right)_{ij} = O_P(\{\log p_n/n\}^{1/2}),$$

which completes the proof of the lemma. □

LEMMA 3. *Let $\mathbf{S}$ be a sample covariance matrix of a random sample $\mathbf{y}_{i 1\leq i\leq n}$ with $\mathbf{y}_i \sim N(\mathbf{0}, \boldsymbol{\Sigma}_0)$. Assume $p_n/n \to y \in [0,1)$, $\boldsymbol{\Sigma}_0$ has eigenvalues uniformly bounded as $n \to \infty$, and $\mathbf{A} = \mathbf{A}_0 + \Delta_1$, $\mathbf{B} = \mathbf{B}_0 + \Delta_2$ are such that the constant matrices $\|\mathbf{A}_0\|, \|\mathbf{B}_0\| = O(1)$, with $\|\Delta_1\|, \|\Delta_2\| = o_P(1)$. Then we still have $\max_{i,j}|(\mathbf{A}(\mathbf{S} - \boldsymbol{\Sigma}_0)\mathbf{B})_{ij}| = O_P(\{\log p_n/n\}^{1/2})$.*

PROOF. Consider

(5.3) $$\mathbf{A}(\mathbf{S} - \boldsymbol{\Sigma}_0)\mathbf{B} = K_1 + K_2 + K_3 + K_4,$$

where $K_1 = \mathbf{A}_0(\mathbf{S} - \boldsymbol{\Sigma}_0)\mathbf{B}_0$, $K_2 = \Delta_1(\mathbf{S} - \boldsymbol{\Sigma}_0)\mathbf{B}_0$, $K_3 = \mathbf{A}_0(\mathbf{S} - \boldsymbol{\Sigma}_0)\Delta_2$ and $K_4 = \Delta_1(\mathbf{S} - \boldsymbol{\Sigma}_0)\Delta_2$. Now, $\max_{i,j}|(K_1)_{ij}| = O_P(\{\log p_n/n\}^{1/2})$ by Lemma 2. Consider $K_2$. Suppose the maximum element of the matrix is at the $(i,j)$th position. Consider $((\mathbf{S} - \boldsymbol{\Sigma}_0)\mathbf{B}_0)_{ij}$, the $(i,j)$th element of $(\mathbf{S} - \boldsymbol{\Sigma}_0)\mathbf{B}_0$. Since each element in $\mathbf{S} - \boldsymbol{\Sigma}_0$ has a rate $O_P(n^{-1/2})$, the $i$th row of $\mathbf{S} - \boldsymbol{\Sigma}_0$ has a norm of $O_P(\{p_n/n\}^{1/2})$. Also, the $j$th column of $\mathbf{B}_0$ has $\|\mathbf{B}_0\mathbf{e}_j\| \leq \|\mathbf{B}_0\| = O(1)$. Hence, $((\mathbf{S} - \boldsymbol{\Sigma}_0)\mathbf{B}_0)_{ij} = O_P(\{p_n/n\}^{1/2})$.

Hence, we can find $c_n = o(\{n/p_n\}^{1/2})$ such that each element in $c_n\mathbf{B}_0^T(\mathbf{S} - \boldsymbol{\Sigma}_0)$ has an order larger than that in $\Delta_1$, since $\|\Delta_1\| = o_P(1)$ implies that each element in $\Delta_1$ is also $o_P(1)$ by Lemma 1.

Then suitable choice of $c_n$ leads to

(5.4) $$\max_{i,j}|(\Delta_1(\mathbf{S} - \boldsymbol{\Sigma}_0)\mathbf{B}_0)_{ij}| \leq c_n \max_k |(\mathbf{B}_0^T(\mathbf{S} - \boldsymbol{\Sigma}_0)^2\mathbf{B}_0)_{kk}|.$$

At the same time, Theorem 5.10 in Bai and Silverstein (2006) implies that, for $\mathbf{y}_i \sim N(\mathbf{0}, \boldsymbol{\Sigma}_0)$ and $p_n/n \to y \in (0,1)$, with probability one,

$$-2\sqrt{y} - y \leq \liminf_{n\to\infty} \lambda_{\min}(\boldsymbol{\Sigma}_0^{-1/2}\mathbf{S}\boldsymbol{\Sigma}_0^{-1/2} - \mathbf{I})$$
$$\leq \limsup_{n\to\infty} \lambda_{\max}(\boldsymbol{\Sigma}_0^{-1/2}\mathbf{S}\boldsymbol{\Sigma}_0^{-1/2} - \mathbf{I}) \leq 2\sqrt{y} + y.$$

Hence, if we have $p_n/n = o(1)$, we must have $\|\boldsymbol{\Sigma}_0^{-1/2}\mathbf{S}\boldsymbol{\Sigma}_0^{-1/2} - \mathbf{I}\| = o_P(1)$, or it will contradict the above. It means that $\|\mathbf{S} - \boldsymbol{\Sigma}_0\| = o_P(1)$ since $\boldsymbol{\Sigma}_0$

has eigenvalues uniformly bounded. Or, if $p_n/n \to y \in (0,1)$, then we have $\|\mathbf{S} - \boldsymbol{\Sigma}_0\| = O_P(1)$ by the above.

Since $\mathbf{S} - \boldsymbol{\Sigma}_0$ is symmetric, we can find a rotation matrix $\mathbf{Q}$ (i.e., $\mathbf{Q}^T\mathbf{Q} = \mathbf{Q}\mathbf{Q}^T = I$) so that

$$\mathbf{S} - \boldsymbol{\Sigma}_0 = \mathbf{Q}\Lambda\mathbf{Q}^T,$$

where $\Lambda$ is a diagonal matrix with real entries. Then we are free to control $c_n$ again so as to satisfy further that $c_n\|\Lambda\|^2 = o_P(\|\Lambda\|)$, since $\|\Lambda\| = \|\mathbf{S} - \boldsymbol{\Sigma}_0\| = O_P(1)$ at most. Hence,

$$\begin{aligned}
c_n \max_k |(\mathbf{B}_0^T (\mathbf{S} - \boldsymbol{\Sigma}_0)^2 \mathbf{B}_0)_{kk}| &= \max_k |(\mathbf{B}_0^T \mathbf{Q} c_n \Lambda^2 \mathbf{Q}^T \mathbf{B}_0)_{kk}| \\
&\leq \max_k |(\mathbf{B}_0^T \mathbf{Q} \Lambda \mathbf{Q}^T \mathbf{B}_0)_{kk}| \\
&= \max_k |(\mathbf{B}_0^T (\mathbf{S} - \boldsymbol{\Sigma}_0) \mathbf{B}_0)_{kk}| \\
&= O_P(\{\log p_n/n\}^{1/2}),
\end{aligned}$$

where the last line used the previous proof for constant matrix $\mathbf{B}_0$. Hence, combining this with (5.4), we have $\max_{i,j}|(K_2)_{ij}| = O_P(\{\log p_n/n\}^{1/2})$. Similar arguments go for $K_3$ and $K_4$. $\square$

PROOF OF THEOREM 1. The main idea of the proof is inspired by Fan and Li (2001) and Rothman et al. (2008). Let $U$ be a symmetric matrix of size $p_n$, $\mathbf{D}_U$ be its diagonal matrix and $\mathbf{R}_U = U - \mathbf{D}_U$ be its off-diagonal matrix. Set $\Delta_U = \alpha_n \mathbf{R}_U + \beta_n \mathbf{D}_U$. We would like to show that, for $\alpha_n = (s_{n1} \log p_n/n)^{1/2}$ and $\beta_n = (p_n \log p_n/n)^{1/2}$, and for a set $\mathcal{A}$ defined as $\mathcal{A} = \{U : \|\Delta_U\|_F^2 = C_1^2 \alpha_n^2 + C_2^2 \beta_n^2\}$,

$$P\Big(\inf_{U \in \mathcal{A}} q_1(\boldsymbol{\Omega}_0 + \Delta_U) > q_1(\boldsymbol{\Omega}_0)\Big) \to 1,$$

for sufficiently large constants $C_1$ and $C_2$. This implies that there is a local minimizer in $\{\boldsymbol{\Omega}_0 + \Delta_U : \|\Delta_U\|_F^2 \leq C_1^2 \alpha_n^2 + C_2^2 \beta_n^2\}$ such that $\|\hat{\boldsymbol{\Omega}} - \boldsymbol{\Omega}_0\|_F = O_P(\alpha_n + \beta_n)$ for sufficiently large $n$, since $\Omega_0 + \Delta_U$ is positive definite. This is shown by noting that

$$\lambda_{\min}(\Omega_0 + \Delta_U) \geq \lambda_{\min}(\Omega_0) + \lambda_{\min}(\Delta_U) \geq \lambda_{\min}(\Omega_0) - \|\Delta_U\|_F > 0,$$

since $\Omega_0$ has eigenvalues uniformly bounded away from 0 and $\infty$ by condition (A), and $\|\Delta_U\|_F = O(\alpha_n + \beta_n) = o(1)$.

Consider, for $\boldsymbol{\Sigma} = \boldsymbol{\Sigma}_0 + \Delta_U$, the difference

$$q_1(\boldsymbol{\Omega}) - q_1(\boldsymbol{\Omega}_0) = I_1 + I_2 + I_3,$$

where

$$I_1 = \text{tr}(\mathbf{S}\mathbf{\Omega}) - \log|\mathbf{\Omega}| - (\text{tr}(\mathbf{S}\mathbf{\Omega}_0) - \log|\mathbf{\Omega}_0|),$$

$$I_2 = \sum_{(i,j)\in S_1^c} (p_{\lambda_{n1}}(|\omega_{ij}|) - p_{\lambda_{n1}}(|\omega_{ij}^0|)),$$

$$I_3 = \sum_{(i,j)\in S_1, i\neq j} (p_{\lambda_{n1}}(|\omega_{ij}|) - p_{\lambda_{n1}}(|\omega_{ij}^0|)).$$

It is sufficient to show that the difference is positive asymptotically with probability tending to 1. Using Taylor's expansion with the integral remainder, we have $I_1 = K_1 + K_2$, where

(5.5)
$$K_1 = \text{tr}((\mathbf{S} - \mathbf{\Sigma}_0)\Delta_U),$$

$$K_2 = \text{vec}(\Delta_U)^T \left\{ \int_0^1 g(v, \mathbf{\Omega}_v)(1-v)\,dv \right\} \text{vec}(\Delta_U)$$

with the definitions $\mathbf{\Omega}_v = \mathbf{\Omega}_0 + v\Delta_U$, and $g(v, \mathbf{\Omega}_v) = \mathbf{\Omega}_v^{-1} \otimes \mathbf{\Omega}_v^{-1}$. Now,

$$K_2 \geq \int_0^1 (1-v)\min_{0\leq v\leq 1} \lambda_{\min}(\mathbf{\Omega}_v^{-1} \otimes \mathbf{\Omega}_v^{-1})\,dv \cdot \|\text{vec}(\Delta_U)\|^2$$

$$= \|\text{vec}(\Delta_U)\|^2/2 \cdot \min_{0\leq v\leq 1} \lambda_{\max}^{-2}(\mathbf{\Omega}_v)$$

$$\geq \|\text{vec}(\Delta_U)\|^2/2 \cdot (\|\mathbf{\Omega}_0\| + \|\Delta_U\|)^{-2}$$

$$\geq (C_1^2\alpha_n^2 + C_2^2\beta_n^2)/2 \cdot (\tau_1^{-1} + o(1))^{-2},$$

where we used $\|\Delta_U\| \leq C_1\alpha_n + C_2\beta_n = O((\log p_n)^{(1-k)/2}) = o(1)$ by our assumption.

Consider $K_1$. It is clear that $|K_1| \leq L_1 + L_2$, where

$$L_1 = \left| \sum_{(i,j)\in S_1} (\mathbf{S} - \mathbf{\Sigma}_0)_{ij}(\Delta_U)_{ij} \right|,$$

$$L_2 = \left| \sum_{(i,j)\in S_1^c} (\mathbf{S} - \mathbf{\Sigma}_0)_{ij}(\Delta_U)_{ij} \right|.$$

Using Lemmas 1 and 2, we have

$$L_1 \leq (s_{n1} + p_n)^{1/2} \max_{i,j} |(\mathbf{S} - \mathbf{\Sigma}_0)_{ij}| \cdot \|\Delta_U\|_F$$

$$\leq O_P(\alpha_n + \beta_n) \cdot \|\Delta_U\|_F$$

$$= O_P(C_1\alpha_n^2 + C_2\beta_n^2).$$

This is dominated by $K_2$ when $C_1$ and $C_2$ are sufficiently large.

Now, consider $I_2 - L_2$ for penalties other than $L_1$. Since $\|\Delta_U\|_F^2 = C_1^2 \alpha_n^2 + C_2^2 \beta_n^2$ on $\mathcal{A}$, we have that $|\omega_{ij}| = O(C_1 \alpha_n + C_2 \beta_n) = o(1)$ for all $(i,j) \in S_1^c$. Also, note that the condition on $\lambda_{n1}$ ensures that, for $(i,j) \in S_1^c$, $|\omega_{ij}| = O(\alpha_n + \beta_n) = o(\lambda_{n1})$. Hence, by condition (C), for all $(i,j) \in S_1^c$, we can find a constant $k_1 > 0$ such that

$$p_{\lambda_{n1}}(|\omega_{ij}|) \geq \lambda_{n1} k_1 |\omega_{ij}|.$$

This implies that

$$I_2 = \sum_{(i,j) \in S_1^c} p_{\lambda_{n1}}(|\omega_{ij}|) \geq \lambda_{n1} k_1 \sum_{(i,j) \in S_1^c} |\omega_{ij}|.$$

Hence,

$$\begin{aligned} I_2 - L_2 &\geq \sum_{(i,j) \in S_1^c} \{\lambda_{n1} k_1 |\omega_{ij}| - |(\mathbf{S} - \boldsymbol{\Sigma}_0)_{ij}| \cdot |\omega_{ij}|\} \\ &\geq \sum_{(i,j) \in S_1^c} [\lambda_{n1} k_1 - O_P(\{\log p_n/n\}^{1/2})] \cdot |\omega_{ij}| \\ &= \lambda_{n1} \sum_{(i,j) \in S_1^c} [k_1 - O_P(\lambda_{n1}^{-1} \{\log p_n/n\}^{1/2})] \cdot |\omega_{ij}|. \end{aligned}$$

With the assumption that $(p_n + s_{n1}) \log p_n / n = O(\lambda_{n1}^2)$, we see from the above that $I_2 - L_2 \geq 0$ since $O_P(\lambda_{n1}^{-1} \{\log p_n/n\}^{1/2}) = o_P(1)$, using $\log p_n/n = o((p_n + s_{n1}) \log p_n/n) = o(\lambda_{n1}^2)$.

For the $L_1$-penalty, since we have $\max_{i \neq j} |\mathbf{S} - \boldsymbol{\Sigma}_0| = O_P((\log p_n/n)^{1/2})$ by Lemma 2, we can find a positive $W = O_P(1)$ such that

$$\max_{i \neq j} |\mathbf{S} - \boldsymbol{\Sigma}_0| = W(\log p_n/n)^{1/2}.$$

Then we can set $\lambda_{n1} = 2W(\log p_n/n)^{1/2}$ or one with order greater than $(\log p_n/n)^{1/2}$, and the above arguments are still valid, so that $I_2 - L_2 > 0$.

Now, with $L_1$ dominated by $K_2$ and $I_2 - L_2 \geq 0$, the proof completes if we can show that $I_3$ is also dominated by $K_2$, since we have proved that $K_2 > 0$. Using Taylor's expansion, we can arrive at

$$|I_3| \leq \min(C_1, C_2)^{-1} \cdot O(1) \cdot (C_1^2 \alpha_n^2 + C_2^2 \beta_n^2) + o(1) \cdot (C_1^2 \alpha_n^2 + C_2^2 \beta_n^2),$$

where $o(1)$ and $O(1)$ are the terms independent of $C_1$ and $C_2$. By condition (B), we have

$$|I_3| = C \cdot O(\alpha_n^2 + \beta_n^2) + C^2 \cdot o(\alpha_n^2 + \beta_n^2),$$

which is dominated by $K_2$ with large enough constants $C_1$ and $C_2$. This completes the proof of the theorem. $\square$

PROOF OF THEOREM 2. For $\mathbf{\Omega}$ a minimizer of (1.1), the derivative for $q_1(\mathbf{\Omega})$ w.r.t. $\omega_{ij}$ for $(i,j) \in S_2^c$ is

$$\frac{\partial q_1(\mathbf{\Omega})}{\partial \omega_{ij}} = 2(s_{ij} - \sigma_{ij} + p'_{\lambda_{n1}}(|\omega_{ij}|)\operatorname{sgn}(\omega_{ij})),$$

where $\operatorname{sgn}(a)$ denotes the sign of $a$. If we can show that the sign of $\partial q_1(\mathbf{\Omega})/\partial \omega_{ij}$ depends on $\operatorname{sgn}(\omega_{ij})$ only with probability tending to 1, the optimum will be at 0, so that $\hat{\omega}_{ij} = 0$ for all $(i,j) \in S_2^c$ with probability tending to 1. We need to estimate the order of $s_{ij} - \sigma_{ij}$ independent of $i$ and $j$.

Decompose $s_{ij} - \sigma_{ij} = I_1 + I_2$, where

$$I_1 = s_{ij} - \sigma_{ij}^0, \qquad I_2 = \sigma_{ij}^0 - \sigma_{ij}.$$

By Lemma 2 or Lemma A.3 of Bickel and Levina (2008b), it follows that $\max_{i,j}|I_1| = O_P(\{\log p_n/n\}^{1/2})$. It remains to estimate the order of $I_2$.

By Lemma 1, $|\sigma_{ij} - \sigma_{ij}^0| \leq \|\mathbf{\Sigma} - \mathbf{\Sigma}_0\|$, which has order

$$\|\mathbf{\Sigma} - \mathbf{\Sigma}_0\| = \|\mathbf{\Sigma}(\mathbf{\Omega} - \mathbf{\Omega}_0)\mathbf{\Sigma}_0\| \leq \|\mathbf{\Sigma}\| \cdot \|\mathbf{\Omega} - \mathbf{\Omega}_0\| \cdot \|\mathbf{\Sigma}_0\| = O(\|\mathbf{\Omega} - \mathbf{\Omega}_0\|),$$

where we used condition (A) to get $\|\mathbf{\Sigma}_0\| = O(1)$, and using $\eta_n \to 0$ so that $\lambda_{\min}(\mathbf{\Omega} - \mathbf{\Omega}_0) = o(1)$ for $\|\mathbf{\Omega} - \mathbf{\Omega}_0\| = O(\eta_n^{1/2})$,

$$\|\mathbf{\Sigma}\| = \lambda_{\min}^{-1}(\mathbf{\Omega}) \leq (\lambda_{\min}(\mathbf{\Omega}_0) + \lambda_{\min}(\mathbf{\Omega} - \mathbf{\Omega}_0))^{-1}$$
$$= (O(1) + o(1))^{-1} = O(1).$$

Hence, $\|\mathbf{\Omega} - \mathbf{\Omega}_0\| = O(\eta_n^{1/2})$ implies $|I_2| = O(\eta_n^{1/2})$.

Combining the last two results yields that

$$\max_{i,j}|s_{ij} - \sigma_{ij}| = O_P(|s_{ij} - \sigma_{ij}^0| + \eta_n^{1/2})$$
$$= O_P(\{\log p_n/n\}^{1/2} + \eta_n^{1/2}).$$

By conditions (C) and (D), we have

$$p'_{\lambda_{n1}}(|\omega_{ij}|) = C_3 \lambda_{n1}$$

for $\omega_{ij}$ in a small neighborhood of 0 (excluding 0 itself) and some positive constant $C_3$. Hence, if $\omega_{ij}$ lies in a small neighborhood of 0, we need to have $\log p_n/n + \eta_n = O(\lambda_{n1}^2)$ in order to have the sign of $\partial q_1(\mathbf{\Omega})/\partial \omega_{ij}$ depends on $\operatorname{sgn}(\omega_{ij})$ only with probability tending to 1. The proof of the theorem is complet. $\square$

PROOF OF THEOREM 3. Because of the similarity between (2.4) and (1.1), the Frobenius norm result has nearly identical proof as Theorem 1, except that we now set $\Delta_U = \alpha_n U$. For the operator norm result, we refer readers to the proof of Theorem 2 of Rothman et al. (2008). $\square$

PROOF OF THEOREM 5. The proof is similar to that of Theorem 1. We only sketch briefly the proof, pointing out the important differences.

Let $\alpha_n = (s_{n2} \log p_n/n)^{1/2}$ and $\beta_n = (p_n \log p_n/n)^{1/2}$, and define $\mathcal{A} = \{U : \|\Delta_U\|_F^2 = C_1^2 \alpha_n^2 + C_2^2 \beta_n^2\}$. Want to show

$$P\Big(\inf_{U \in \mathcal{A}} q_2(\boldsymbol{\Sigma}_0 + \Delta_U) > q_2(\boldsymbol{\Sigma}_0)\Big) \to 1$$

for sufficiently large constants $C_1$ and $C_2$.

For $\boldsymbol{\Sigma} = \boldsymbol{\Sigma}_0 + \Delta_U$, the difference

$$q_2(\boldsymbol{\Sigma}) - q_2(\boldsymbol{\Sigma}_0) = I_1 + I_2 + I_3,$$

where

$$I_1 = \operatorname{tr}(S\boldsymbol{\Omega}) + \log|\boldsymbol{\Sigma}| - (\operatorname{tr}(S\boldsymbol{\Omega}_0) + \log|\boldsymbol{\Sigma}_0|),$$

$$I_2 = \sum_{(i,j) \in S_2^c} (p_{\lambda_{n2}}(|\sigma_{ij}|) - p_{\lambda_{n2}}(|\sigma_{ij}^0|)),$$

$$I_3 = \sum_{(i,j) \in S_2, i \neq j} (p_{\lambda_{n2}}(|\sigma_{ij}|) - p_{\lambda_{n2}}(|\sigma_{ij}^0|))$$

with $I_1 = K_1 + K_2$, where

(5.6)
$$K_1 = -\operatorname{tr}((\mathbf{S} - \boldsymbol{\Sigma}_0)\boldsymbol{\Omega}_0 \Delta_U \boldsymbol{\Omega}_0) = -\operatorname{tr}((\mathbf{S}_{\boldsymbol{\Omega}_0} - \boldsymbol{\Omega}_0)\Delta_U),$$

$$K_2 = \operatorname{vec}(\Delta_U)^T \left\{\int_0^1 g(v, \boldsymbol{\Sigma}_v)(1-v)\,dv\right\} \operatorname{vec}(\Delta_U),$$

and $\boldsymbol{\Sigma}_v = \boldsymbol{\Sigma}_0 + v\Delta_U$, $\mathbf{S}_{\boldsymbol{\Omega}_0}$ is the sample covariance matrix of a random sample $\{\mathbf{x}_i\}_{1 \leq i \leq n}$ having $\mathbf{x}_i \sim N(\mathbf{0}, \boldsymbol{\Omega}_0)$. Also,

(5.7) $\quad g(v, \boldsymbol{\Sigma}_v) = \boldsymbol{\Sigma}_v^{-1} \otimes \boldsymbol{\Sigma}_v^{-1} S \boldsymbol{\Sigma}_v^{-1} + \boldsymbol{\Sigma}_v^{-1} S \boldsymbol{\Sigma}_v^{-1} \otimes \boldsymbol{\Sigma}_v^{-1} - \boldsymbol{\Sigma}_v^{-1} \otimes \boldsymbol{\Sigma}_v^{-1}.$

The treatment of $K_2$ is different from that in Theorem 1. By condition (A), and $(p_n + s_{n2})(\log p_n)^k/n = O(1)$ for some $k > 1$, we have

$$\|v\Delta_U \Omega_0\| \leq \|\Delta_U\| \|\boldsymbol{\Omega}_0\| \leq \tau_1^{-1}(C_1 \alpha_n + C_2 \beta_n) = O((\log p_n)^{1-k}) = o(1).$$

Thus, we can use the Neumann series expansion to arrive at

$$\boldsymbol{\Sigma}_v^{-1} = \boldsymbol{\Omega}_0(I + v\Delta_U \boldsymbol{\Omega}_0)^{-1} = \boldsymbol{\Omega}_0(I - v\Delta_U \boldsymbol{\Omega}_0 + o(1)),$$

where the little $o$ (or $o_P$, $O$ or $O_P$ in any matrix expansions in the remainder of this proof) represents a function of the $L_2$ norm of the residual matrix in the expansion. That is, $\boldsymbol{\Sigma}_v^{-1} = \boldsymbol{\Omega}_0 + O_P(\alpha_n + \beta_n)$, and $\|\boldsymbol{\Sigma}_v^{-1}\| = \tau_1^{-1} + O_P(\alpha_n + \beta_n)$. With $\mathbf{S}_I$ defined as the sample covariance matrix formed from a random sample $\{\mathbf{x}_i\}_{1 \leq i \leq n}$ having $\mathbf{x}_i \sim N(\mathbf{0}, I)$,

$$\|\mathbf{S} - \boldsymbol{\Sigma}_0\| = O_P(\|\mathbf{S}_I - I\|) = o_P(1)$$

(see arguments in Lemma 3). These entail

$$\mathbf{S}\mathbf{\Sigma}_v^{-1} = (\mathbf{S} - \mathbf{\Sigma}_0)\mathbf{\Sigma}_v^{-1} + \mathbf{\Sigma}_0\mathbf{\Sigma}_v^{-1} = o_P(1) + I + O_P(\alpha_n + \beta_n) = I + o_P(1).$$

Combining these results, we have

$$g(v, \mathbf{\Sigma}_v) = \mathbf{\Omega}_0 \otimes \mathbf{\Omega}_0 + O_P(\alpha_n + \beta_n).$$

Consequently,

$$K_2 = \text{vec}(\Delta_U)^T \left\{ \int_0^1 \mathbf{\Omega}_0 \otimes \mathbf{\Omega}_0 (1 + o_P(1))(1 - v)\, dv \right\} \text{vec}(\Delta_U)$$
$$\geq \lambda_{\min}(\mathbf{\Omega}_0 \otimes \mathbf{\Omega}_0) \| \text{vec}(\Delta_U) \|^2 / 2 \cdot (1 + o_P(1))$$
$$= \tau_1^{-2}(C_1^2 \alpha_n^2 + C_2^2 \beta_n^2)/2 \cdot (1 + o_P(1)).$$

All other terms are dealt with similarly as in the proof of Theorem 1, and hence we omit them. □

PROOF OF THEOREM 6. The proof is similar to that of Theorem 2. We only show the main differences.

It is easy to show

$$\frac{\partial q_2(\mathbf{\Sigma})}{\partial \sigma_{ij}} = 2(-(\mathbf{\Omega}\mathbf{S}\mathbf{\Omega})_{ij} + \omega_{ij} + p'_{\lambda_n}(|\sigma_{ij}|)\text{sgn}(\sigma_{ij})).$$

Our aim is to estimate the order of $|(-\mathbf{\Omega}\mathbf{S}\mathbf{\Omega} + \mathbf{\Omega})_{ij}|$, finding an upper bound which is independent of both $i$ and $j$.

Write

$$-\mathbf{\Omega}\mathbf{S}\mathbf{\Omega} + \mathbf{\Omega} = I_1 + I_2,$$

where $I_1 = -\mathbf{\Omega}(\mathbf{S} - \mathbf{\Sigma}_0)\mathbf{\Omega}$ and $I_2 = \mathbf{\Omega}(\mathbf{\Sigma} - \mathbf{\Sigma}_0)\mathbf{\Omega}$. Since

$$\|\mathbf{\Omega}\| = \lambda_{\min}^{-1}(\mathbf{\Sigma}) \leq (\lambda_{\min}(\mathbf{\Sigma}_0) + \lambda_{\min}(\mathbf{\Sigma} - \mathbf{\Sigma}_0))^{-1}$$
$$= \tau_1^{-1} + o(1),$$

we have

$$\mathbf{\Omega} = \mathbf{\Omega}_0 + (\mathbf{\Omega} - \mathbf{\Omega}_0) = \mathbf{\Omega}_0 - \mathbf{\Omega}(\mathbf{\Sigma} - \mathbf{\Sigma}_0)\mathbf{\Omega}_0 = \mathbf{\Omega}_0 + \Delta,$$

where $\|\Delta\| \leq \|\mathbf{\Omega}\| \cdot \|\mathbf{\Sigma} - \mathbf{\Sigma}_0\| \cdot \|\mathbf{\Omega}_0\| = O(\eta_n^{1/2}) = o(1)$ by Lemma 1, with $\|\mathbf{\Sigma} - \mathbf{\Sigma}_0\|^2 = O(\eta_n)$. Hence, we can apply Lemma 3 and conclude that $\max_{i,j} |(I_1)_{ij}| = O_P(\{\log p_n/n\}^{1/2})$.

For $I_2$, we have

$$\max_{i,j} |(I_2)_{ij}| \leq \|\mathbf{\Omega}\| \cdot \|\mathbf{\Sigma} - \mathbf{\Sigma}_0\| \cdot \|\mathbf{\Omega}\| = O(\|\mathbf{\Sigma} - \mathbf{\Sigma}_0\|) = O(\eta_n^{1/2}).$$

Hence, we have

$$\max_{i,j}|(-\mathbf{\Omega S\Omega}+\mathbf{\Omega})_{ij}|=O(\{\log p_n/n\}^{1/2}+\eta_n^{1/2}).$$

The rest goes similar to the proof of Theorem 2, and is omitted. □

PROOF OF THEOREM 7. The proof is nearly identical to that of Theorem 5, except that we now set $\Delta_U = \alpha_n U$. The fact that $(\hat{\mathbf{\Gamma}}_\mathbf{S})_{ii} = 1 = \gamma_{ii}^0$ has no estimation error eliminates an order $(p_n \log p_n/n)^{1/2}$ that contributes from estimating $\text{tr}((\hat{\mathbf{\Gamma}}_\mathbf{S} - \mathbf{\Gamma}_0)\mathbf{\Psi}_0 \Delta_U \mathbf{\Psi}_0)$ for (3.2). This is why we can estimate a sparse correlation matrix more accurately.

For the operator norm result, we refer readers to the proof of Theorem 2 of Rothman et al. (2008). □

PROOF OF THEOREM 10. For $(\mathbf{T}, \mathbf{D})$ a minimizer of (4.2), the derivative for $q_3(\mathbf{T}, \mathbf{D})$ w.r.t. $t_{ij}$ for $(i,j) \in S_3^c$ is

$$\frac{\partial q_3(\mathbf{T}, \mathbf{D})}{\partial t_{ij}} = 2((\mathbf{ST}^T\mathbf{D}^{-1})_{ji} + p'_{\lambda_{n3}}(|t_{ij}|)\text{sgn}(t_{ij})).$$

Now $\mathbf{ST}^T\mathbf{D}^{-1} = I_1 + I_2 + I_3 + I_4$, where

$$I_1 = (\mathbf{S}-\mathbf{\Sigma}_0)\mathbf{T}^T\mathbf{D}^{-1}, \qquad I_2 = \mathbf{\Sigma}_0(\mathbf{T}-\mathbf{T}_0)^T\mathbf{D}^{-1},$$
$$I_3 = \mathbf{\Sigma}_0\mathbf{T}_0^T(\mathbf{D}^{-1}-\mathbf{D}_0^{-1}), \qquad I_4 = \mathbf{\Sigma}_0\mathbf{T}_0^T\mathbf{D}_0^{-1}.$$

By the MCD (4.1), $I_4 = \mathbf{T}_0^{-1}$. Since $i > j$ for $(i,j) \in S_3^c$, we must have $(\mathbf{T}_0^{-1})_{ji} = 0$. Hence, we can ignore $I_4$.

Since $\|\mathbf{T}-\mathbf{T}_0\|^2 = O(\eta_n)$ and $\|\mathbf{D}-\mathbf{D}_0\|^2 = O(\zeta_n)$ with $\eta_n, \zeta_n = o(1)$, and by condition (A) we can easily show $\|\mathbf{D}^{-1}-\mathbf{D}_0^{-1}\| = O(\|\mathbf{D}-\mathbf{D}_0\|) = O(\zeta_n^{1/2})$. Then we can apply Lemma 3 to show that $\max_{ij}|(I_1)_{ij}| = (\log p_n/n)^{1/2}$.

For $I_2$, we have $\max_{ij}|(I_2)_{ij}| \leq \|\mathbf{\Sigma}_0\| \cdot \|\mathbf{T}-\mathbf{T}_0\| \cdot \|\mathbf{D}^{-1}\| = O(\eta_n^{1/2})$. And finally, $\max_{ij}|(I_3)_{ij}| \leq \|\mathbf{\Sigma}_0\| \cdot \|\mathbf{T}_0\| \cdot \|\mathbf{D}^{-1}-\mathbf{D}_0^{-1}\| = O(\zeta_n^{1/2})$.

With all these, we have $\max_{(i,j)\in S_3^c}|(\mathbf{ST}^T\mathbf{D}^{-1})_{ji}|^2 = \log p_n/n + \eta_n + \zeta_n$. The rest of the proof goes like that of Theorems 2 or 6. □

## REFERENCES

BAI, Z. and SILVERSTEIN, J. W. (2006). *Spectral Analysis of Large Dimensional Random Matrices*. Science Press, Beijing.

BICKEL, P. J. and LEVINA, E. (2008a). Covariance regularization by thresholding. *Ann. Statist.* **36** 2577–2604. MR2485008

BICKEL, P. J. and LEVINA, E. (2008b). Regularized estimation of large covariance matrices. *Ann. Statist.* **36** 199–227. MR2387969

CAI, T., ZHANG, C.-H. and ZHOU, H. (2008). Optimal rates of convergence for covariance matrix estimation. Technical report, The Wharton School, Univ. Pennsylvania.

D'ASPREMONT, A., BANERJEE, O. and EL GHAOUI, L. (2008). First-order methods for sparse covariance selection. *SIAM J. Matrix Anal. Appl.* **30** 56–66. MR2399568

DEMPSTER, A. P. (1972). Covariance selection. *Biometrics* **28** 157–175.

DIGGLE, P. and VERBYLA, A. (1998). Nonparametric estimation of covariance structure in longitudinal data. *Biometrics* **54** 401–415.

EL KAROUI, N. (2008). Operator norm consistent estimation of a large dimensional sparse covariance matrices. *Ann. Statist.* **36** 2717–2756. MR2485011

FAN, J., FENG, Y. and WU, Y. (2009). Network exploration via the adaptive LASSO and SCAD penalties. *Ann. Appl. Stat.* **3** 521–541.

FAN, J. and LI, R. (2001). Variable selection via nonconcave penalized likelihood and its oracle properties. *J. Amer. Statist. Assoc.* **96** 1348–1360. MR1946581

FAN, J. and PENG, H. (2004). Nonconcave penalized likelihood with a diverging number of parameters. *Ann. Statist.* **32** 928–961. MR2065194

FRIEDMAN, J., HASTIE, T. and TIBSHIRANI, R. (2008). Sparse inverse covariance estimation with the graphical LASSO. *Biostatistics* **9** 432–441.

HUANG, J., HOROWITZ, J. and MA, S. (2008). Asymptotic properties of bridge estimators in sparse high-dimensional regression models. *Ann. Statist.* **36** 587–613. MR2396808

HUANG, J., LIU, N., POURAHMADI, M. and LIU, L. (2006). Covariance matrix selection and estimation via penalised normal likelihood. *Biometrika* **93** 85–98. MR2277742

LEVINA, E., ROTHMAN, A. J. and ZHU, J. (2008). Sparse estimation of large covariance matrices via a nested Lasso penalty. *Ann. Appl. Stat.* **2** 245–263.

MEIER, L., VAN DE GEER, S. and BÜHLMANN, P. (2008). The group Lasso for logistic regression. *J. R. Stat. Soc. Ser. B Stat. Methodol.* **70** 53–71. MR2412631

MEINSHAUSEN, N. and BÜHLMANN, P. (2006). High-dimensional graphs and variable selection with the Lasso. *Ann. Statist.* **34** 1436–1462. MR2278363

POURAHMADI, M. (1999). Joint mean-covariance models with applications to longitudinal data: Unconstrained parameterisation. *Biometrika* **86** 677–690. MR1723786

RAVIKUMAR, P., LAFFERTY, J., LIU, H. and WASSERMAN, L. (2007). Sparse additive models. In *Advances in Neural Information Processing Systems* **20**. MIT Press, Cambridge, MA.

ROTHMAN, A. J., BICKEL, P. J., LEVINA, E. and ZHU, J. (2008). Sparse permutation invariant covariance estimation. *Electron. J. Stat.* **2** 494–515. MR2417391

SMITH, M. and KOHN, R. (2002). Parsimonious covariance matrix estimation for longitudinal data. *J. Amer. Statist. Assoc.* **97** 1141–1153. MR1951266

WAGAMAN, A. S. and LEVINA, E. (2008). Discovering sparse covariance structures with the Isomap. *J. Comput. Graph. Statist.* **18**. To appear.

WONG, F., CARTER, C. and KOHN, R. (2003). Efficient estimation of covariance selection models. *Biometrika* **90** 809–830. MR2024759

WU, W. B. and POURAHMADI, M. (2003). Nonparametric estimation of large covariance matrices of longitudinal data. *Biometrika* **94** 1–17. MR2024760

YUAN, M. and LIN, Y. (2007). Model selection and estimation in the Gaussian graphical model. *Biometrika* **90** 831–844. MR2367824

ZHANG, C. H. (2007). Penalized linear unbiased selection. Technical report 2007-003, The Statistics Dept., Rutgers Univ.

ZHAO, P. and YU, B. (2006). On model selection consistency of Lasso. *J. Mach. Learn. Res.* **7** 2541–2563. MR2274449

ZOU, H. (2006). The adaptive Lasso and its oracle properties. *J. Amer. Statist. Assoc.* **101** 1418–1429. MR2279469

ZOU, H. and LI, R. (2008). One-step sparse estimates in nonconcave penalized likelihood models (with discussion). *Ann. Statist.* **36** 1509–1533. MR2435443

Department of Statistics  
London School of Economics  
  and Political Science  
Houghton Street  
London, WC2A 2AE  
United Kingdom  
E-mail: C.Lam2@lse.ac.uk

Department of Operations Research  
  and Financial Engineering  
Princeton University  
Princeton, New Jersey 08544  
USA  
E-mail: jqfan@princeton.edu